\documentclass[12pt]{article}
\usepackage{graphicx}
\usepackage{graphics}
\usepackage{amsthm}
\usepackage{amssymb, amsmath}
\title{Generalizations of Russell-style integrals} 
\author{Mark W. Coffey\\
Department of Physics\\
Colorado School of Mines\\
Golden, CO  80401\\
USA\\
mcoffey@mines.edu\\}
\date{June 6, 2018}
\begin{document}
\maketitle
\baselineskip=25 pt
\begin{abstract}

First some definite integrals of W. H. L. Russell, almost all with trigonometric function
integrands, are derived, and many generalized.  Then a list is given in Russell-style of
generalizations of integral identities of Amdeberhan and Moll. 
We conclude with a brief and noncomprehensive description of directions for further investigation,
including the significant generalization to elliptic functions.

\end{abstract}
 
\medskip
\baselineskip=15pt
\centerline{\bf Key words and phrases}
\medskip 

definite integrals, Gauss hypergeometric function, Jacobi elliptic function 

\bigskip
\noindent
{\bf 2010 AMS codes}
\newline{26A42, 33C60}

\baselineskip=25pt

\pagebreak
\centerline{\bf 1.  Introduction} 
\medskip

This paper first reviews, evaluates, and extends various definite integrals of W. H. L. Russell
\cite{russell}.  It then generalizes several integrals of \cite{tewodros1} that are stated to be
of ``Russell-style".  Because several of the results are expressed in terms of the Gauss
hypergeometric function ${}_2F_1$ for general argument and rational parameters, they likely
have application in the calculation of hypervolumes and of the terms of Feynman diagrams.
Several directions for further research are mentioned in the last section.

\medskip
{\bf 2. Some integrals of Russell}
\medskip


There are a dozen integrals presented without proofs in the listing \cite{russell}, which we
label successively with a preceding R.
We evaluate the first seven of these integrals, giving generalizations for most of them.  
For several, it is useful to employ the families of integrals contained in the following.
{\newline \bf Lemma 1}.  
There holds for integer $m \geq 0$
$$\int_0^\pi \theta \sin \theta \cos^{2m}\theta d\theta={\pi \over {2m+1}},$$
$$\int_0^\pi \theta \sin^3 \theta \cos^{2m} \theta d\theta={{2\pi} \over {4m^2+8m+3}},$$
$$\int_0^\pi \theta \sin^5 \theta \cos^{2m} \theta d\theta={{8\pi} \over {(2m+1)(2m+3)(2m+5)}},$$
$$\int_0^\pi \theta \sin^7 \theta \cos^{2m} \theta d\theta={{48\pi} \over {(2m+1)(2m+3)(2m+5)(2m+7)}},$$
and
$$\int_0^\pi \theta \sin^{2j+1} \theta \cos^{2m} \theta d\theta={{2^j j! \pi} \over {(2m+1)(2m+3)(2m+5)\cdots(2m+2j+1)}}.$$
There are many relations between these integrals, which may be studied via recursion relations and
otherwise.

Let
$$I_{r,p,q}\equiv \int x^r \sin^p x\cos^q x ~dx.$$
Then the latter entries of Lemma 1 follow from the first as the $r=1$ special case of
$$I_{r,p,q}={1 \over {(p+q)^2}}\left[(p+q)x^r \sin^{p+1}x \cos^{q-1}x+rx^{r-1}\sin^p x \cos^q x \right.$$
$$\left. -r(r-1)I_{r-2,p,q}-rpI_{r-1,p-1,q-1}+(q-1)(p+q)I_{r,p,q-2}\right].$$

Specifically for the first entry of Lemma 1 we have the easy integration by parts, followed by the use of
a Beta function integral,
$$\int_0^\pi \theta \sin \theta \cos^m\theta d\theta={1 \over {m+1}}\int_0^\pi \cos^{m+1} \theta d\theta
-{\pi \over {m+1}}(-1)^{m+1}$$
$$={\sqrt{\pi} \over {2(m+1)}} {{\Gamma \left({m \over 2}+1\right)} \over {\Gamma\left({m \over 2}+{3 \over 2}
\right)}}[1-(-1)^m]-{\pi \over {m+1}}(-1)^{m+1}.$$

The first entry of Lemma 1 is a highly specialized case of the following indefinite integral.
{\newline \bf Lemma 2}.  Up to an additive constant,
$$\int t \sin t \cos^m t ~dt=-{{\cos^{m+1} t} \over {(m+1)(m+2)}}\left[\cos t ~{}_2F_1\left({1 \over 2};
{m \over 2}+1;{m \over 2}+2;\cos^2 t\right) + (m+2)t \right].$$
The form of this result shows that explicit evaluations are possible at least for $t=0$, $\pi/4$, $\pi/2$,
and $\pi$.  This follows since ${}_2F_1(\alpha,\beta;\gamma;0)=1$, ${}_2F_1(\alpha,\beta;\gamma;1)
=\Gamma(\gamma)\Gamma(\gamma-\alpha-\beta)/[\Gamma(\gamma-\alpha)\Gamma(\gamma-\beta)]$ in terms of
the Gamma function, and ${}_2F_1(\alpha,\beta;\gamma;1/2)$ transforms to $2^\alpha {}_2F_1(\alpha,\gamma-\beta;\gamma;-1)$ and, equivalently, to $2^\beta {}_2F_1(\beta,\gamma-\alpha;\gamma;-1)$.
The Gauss hypergeometric function in Lemma 2 has the further property that Clausen's identity applies to ${}_2F_1(\alpha,\beta;\alpha+\beta+1/2;x)$.
One method to develop Lemma 2 is to apply the expansion, for example, for even exponents as we are
interested herein,
$$\cos^{2m} x={1 \over 2^{2m}}\left[\sum_{k=0}^{m-1}2{{2m} \choose k}\cos 2(m-k)x+{{2m} \choose m}\right],$$
wherein $\cos 2(m-k)x=2\cos^2(m-k)x-1$ and ${{2m} \choose m}$ is the central binomial coefficient.

{\it Remark}.  We refrain from doing so here, but it is possible to then consider the analogous integrals
$$\int t \cos t \sin^m t ~dt, ~~~~~~~~ \int t \cos^{2j+1} t \sin^m t ~dt,$$
and analogs of the Russel integrals which are next taken up. 

R1 is given as
$$\int_0^\infty e^{-(r+1)z+xe^{-z}}dz={e^x \over x}\left(1-{r \over x}+r(r-1){1 \over x^2}-\ldots
\right).$$
This integral may be written in terms of the Gamma $\Gamma(x)$ and incomplete $\Gamma(x,y)$
functions \cite{nbs,grad}:
$$\int_0^\infty e^{-(r+1)z+xe^{-z}}dz=\sum_{n=0}^\infty {x^n \over {n!}}\int_0^\infty e^{-(n+r+1)z}dz$$
$$=\sum_{n=0}^\infty {x^n \over {n!}}{1 \over {(n+r+1)}}$$
$$=(-x)^{-(r+1)}[\Gamma(r+1)-\Gamma(r+1,-x)].$$


R2 is:
$$\int_0^\pi {{\theta \sin \theta} \over {1-x^2\cos^2 \theta}}d\theta={\pi \over {2x}}\ln {{(1+x)} \over
{(1-x)}},$$
while R3 is:
$$\int_0^\pi {{\theta \sin \theta} \over {1+x^2\cos^2 \theta}}={\pi \over x}\tan^{-1} x.$$
Hence, using complex variables, these are equivalent, and R2 could just as well be written with
the inverse function $\tanh^{-1}$.  In turn, these integrals are subsumed by the generalization
given for R7 below.  By expanding the denominator of the integrand as geometric series and using 
Lemma 1 we obtain
$$\int_0^\pi {{\theta \sin \theta} \over {1-x^2\cos^2 \theta}}d\theta
=\sum_{n=0}^\infty x^{2n} \int_0^\pi \theta \sin \theta \cos^{2n} \theta d\theta$$
$$=\sum_{n=0}^\infty {x^{2n} \over {2n+1}}={\pi \over x}\tanh^{-1} x.$$

By using other parts of Lemma 1, generalizations of R2/R3 are easily obtained.
For instance:
$$\int_0^\pi {{\theta \sin^3 \theta} \over {1-x^2\cos^2 \theta}}d\theta
=\sum_{n=0}^\infty x^{2n} \int_0^\pi \theta \sin^3 \theta \cos^{2n} \theta d\theta$$
$$=\sum_{n=0}^\infty {x^{2n} \over {(2n+1)(2n+3)}}={\pi \over x^2}\left[1+{{x^2-1} \over x}\tanh^{-1} x
\right]$$
and
$$\int_0^\pi {{\theta \sin^5 \theta} \over {1-x^2\cos^2 \theta}}d\theta
=\sum_{n=0}^\infty x^{2n} \int_0^\pi \theta \sin^5 \theta \cos^{2n} \theta d\theta$$
$$=\sum_{n=0}^\infty {x^{2n} \over {(2n+1)(2n+3)(2n+5)}}={\pi \over x^4}\left[{5 \over 3}x^2-1
+{{(x^2-1)^2} \over x}\tanh^{-1} x\right].$$

R4 is:
$$\int_0^{\pi/2} e^{\cos^3 \theta \cos 3\theta}\cos(\cos^3 \theta \sin 3\theta)d\theta
={\pi \over 2}e^{1/8}.$$
We show the generalization
$$\int_0^{\pi/2} e^{\cos^j \theta \cos j\theta}\cos(\cos^j \theta \sin j\theta)d\theta
={\pi \over 2}e^{1/2^j}.$$
$$\int_0^{\pi/2} e^{\cos^j \theta \cos j\theta}\cos(\cos^j \theta \sin j\theta)d\theta
=\mbox{Re}~ \int_0^{\pi/2} e^{e^{ji\theta}\cos^j \theta}d\theta$$
$$=\mbox{Re}~ \sum_{n=0}^\infty {1 \over {n!}} \int_0^{\pi/2} \cos^{jn} \theta e^{jin\theta}d\theta$$
$$={1 \over 2} \sum_{n=0}^\infty {1 \over {n!}} \int_{-\pi/2}^{\pi/2} \cos^{jn} \theta e^{jin\theta}d\theta$$
$$={\pi \over 2} \sum_{n=0}^\infty {1 \over {n!}}{1 \over 2^{jn}}={\pi \over 2}e^{1/2^j}.$$
The integral in the penultimate line was evaluated by way of the reciprocal Beta function
expression 3.892.2 of \cite{grad} (p.\ 476).

R5 is stated as
$$\int_0^\pi \theta \sin \theta \sqrt{1+x^2 \cos^2 \theta} d\theta
={\pi \over {2x}}[x\sqrt{1+x^2}+\ln(x+\sqrt{1+x^2})].$$

We demonstrate the generalization
$$\int_0^\pi \theta \sin \theta (1+x^2 \cos^2 \theta)^{1/k}d\theta
=\pi {}_2F_1\left({1 \over 2},-{1 \over k};{3 \over 2};-x^2\right),$$
where ${}_2F_1$ is the Gauss hypergeometric function \cite{nbs,grad}.
From binomial expansion within the integrand,
$$\int_0^\pi \theta \sin \theta (1+x^2 \cos^2 \theta)^{1/k}d\theta
=\sum_{j=0}^\infty {{1/k} \choose j}x^{2j}\int_0^\pi \theta \sin \theta \cos^{2j}\theta d\theta$$
$$=\pi \sum_{j=0}^\infty {{1/k} \choose j}{x^{2j} \over {2j+1}}$$
$$=\pi {}_2F_1\left({1 \over 2},-{1 \over k};{3 \over 2};-x^2\right).$$
For the last step we have noted that ${{1/k} \choose j}=(-1)^j (-1/k)_j/j!$ and $1/(2j+1)=(1/2)_j
/(3/2)_j$ in terms of Pochhammer symbol $(a)_j=\Gamma(a+j)/\Gamma(a)$.
Another way to obtain this result is to write
$$\pi \sum_{j=0}^\infty {{1/k} \choose j}{x^{2j} \over {2j+1}}
=\pi \sum_{j=0}^\infty {{1/k} \choose j} x^{2j}\int_0^\infty e^{-(2j+1)t}dt$$
$$=\pi \int_0^\infty e^{-t}(1+e^{-2t}x^2)^{1/k}dt=\pi \int_0^1 (1+x^2 u^2)^{1/k}du$$
$$=\pi {}_2F_1\left({1 \over 2},-{1 \over k};{3 \over 2};-x^2\right).$$

For $k=2$ we have the R5 reduction
$$\int_0^\pi \theta \sin \theta (1+x^2 \cos^2 \theta)^{1/2}d\theta={\pi \over 2}\left[\sqrt{1+x^2}
+{1 \over x}\sinh^{-1}x \right].$$
Here the case 
$${}_2F_1\left({1 \over 2},-{1 \over 2};{3 \over 2};-x^2\right)$$
may be obtained from 
$${}_2F_1\left({1 \over 2},{1 \over 2};{3 \over 2};-x^2\right)
=\sum_{k=0}^\infty (-1)^k {{(2k)!} \over {4^k (k!)^2(2k+1)}}x^{2k} 
={1 \over x}\sinh^{-1}x$$ 
via the relation $(1/2)_j=-2(j-1/2)(-1/2)_j$.  The connection of the middle display expression to
hypergeometric form is by way of the duplication formula for Pochhammer symbols 
$$(\lambda)_n=2^{2n}\left({\lambda \over 2}\right)_n\left({{\lambda+1} \over 2}\right)_n.$$

{\it Remark}.  The integrals of R2, R3, and R5 may be generalized to arbitrary powers of $\sin \theta$
within the integrand as we show below for R6.

Entry R6 is
$$\int_0^\pi \theta \sin \theta [(1+x\cos\theta)^{1/3}+(1-x\cos \theta)^{1/3}]d\theta
={{3\pi} \over {4x}}[(1+x)^{4/3}-(1-x)^{4/3}],$$
whereas we first demonstrate the generalization
$$\int_0^\pi \theta \sin \theta [(1+x\cos\theta)^{1/k}+(1-x\cos \theta)^{1/k}]d\theta
={{\pi k} \over {(k+1)x}}[(1+x)^{(k+1)/k}-(1-x)^{(k+1)/k}].$$
Again binomial expansion and Lemma 1 are applied so that
$$\int_0^\pi \theta \sin \theta [(1+x\cos\theta)^{1/k}+(1-x\cos \theta)^{1/k}]d\theta
=\sum_{j=0}^\infty {{1/k} \choose j} x^j[1+(-1)^j]\int_0^\infty \theta \sin \theta \cos^j \theta d\theta$$
$$=2\sum_{m=0}^\infty {{1/k} \choose {2m}}x^{2m} \int_0^\infty \theta \sin \theta \cos^j \theta d\theta$$
$$=2\pi \sum_{m=0}^\infty {{1/k} \choose {2m}}{x^{2m} \over {(2m+1)}}$$
$$={{\pi k} \over {(k+1)x}}[(1+x)^{(k+1)/k}-(1-x)^{(k+1)/k}].$$
The final summation follows from an easy integration of
$$\sum_{m=0}^\infty {{1/k} \choose {2m}}x^{2m}={1 \over 2}[(1-x)^{1/k}+(1+x)^{1/k}].$$

We have several further generalizations by using Lemma 1, including
$$\int_0^\pi \theta \sin^3 \theta [(1+x\cos\theta)^{1/k}+(1-x\cos \theta)^{1/k}]d\theta
={{4\pi} \over 3} {}_2F_1\left(-{1 \over {2k}},{{k-1} \over {2k}};{5 \over 2};x^2\right).$$
When $k=2$ we have the result
$$\int_0^\pi \theta \sin^3 \theta [(1+x\cos\theta)^{1/2}+(1-x\cos \theta)^{1/2}]d\theta$$
$$={{8\pi} \over {105 x^3}}[\sqrt{1-x}(2+x-8x^2+5x^3)+\sqrt{1+x}(-2+x+8x^2+5x^3)].$$ 
In fact, for $k>2$,
$$\int_0^\pi \theta \sin^3 \theta [(1+x\cos\theta)^{1/k}+(1-x\cos \theta)^{1/k}]d\theta$$
has the form
$${{4\pi} \over 3}{q \over x^3}\left[{{p(x)} \over {(1-x)^{(2k-1)/2}}}+{{\tilde{p}(x)} \over {(1+x)^{(2k-1)/2}}}\right],$$
where $q$ is a rational number, $p$ is a quintic polynomial, and $\tilde{p}(x)=-p(-x)$.  
Moreover, explicitly,
$$p(x)=(x-1)^4 \left(x+{k \over {2k+1)}}\right), 
~~~~ \tilde{p}(x)=(x+1)^4 \left(x-{k \over {2k+1)}}\right).$$

Omitting many other special cases and examples, we present the following general result coming 
from Lemma 1.
\newline{For odd positive integers $p$,}
$$\int_0^\pi \theta \sin^p \theta [(1+x\cos\theta)^{1/k}+(1-x\cos \theta)^{1/k}]d\theta
={{2^{(p+1)/2}[(p-1)/2]!\pi} \over {p!!}} {}_2F_1\left(-{1 \over {2k}},{{k-1} \over {2k}};{{p+2} \over
2};x^2\right).$$

R7 reads
$$\int_0^\pi {{\theta \sin \theta} \over {1+x^4 \cos^4 \theta}}d\theta
={\pi \over {2^{5/2} x}}\ln\left({{1+\sqrt{2}x+x^2} \over {1-\sqrt{2}x+x^2}}\right)+{\pi \over {2^{3/2}x}}
\tan^{-1}{{\sqrt{2}x} \over {1-x^2}},$$
and is also easily generalized.


We have
$$\int_0^\pi {{\theta \sin \theta} \over {1+x^j \cos^j \theta}}d\theta
=\sum_{n=0}^\infty (-1)^j x^{jn}\int_0^\pi \theta \sin \theta \cos^{jn} \theta d\theta$$
$$=\pi \sum_{n=0}^\infty (-1)^j {x^{jn} \over {jn+1}}$$
$$=\pi {{(-1)^j} \over j}\Phi\left(x^j,1,{1 \over j}\right),$$
wherein $\Phi(z,a,b)$ is the Lerch zeta function \cite{nbs,grad}.  For $j=4$ we obtain the R7 reduction 
to $\pi ~{}_2F_1(1/4,1;5/4;x^4)$.

{\it Remarks}.  In regard to R8 one has available the expansions
$$\ln \cos x=\sum_{k=1}^\infty {{2^{2k-1}(2^{2k}-1)B_{2k}} \over {k(2k)!}}x^{2k}$$
$$=-{1 \over 2}\sum_{k=1}^\infty {{\sin^{2k} x} \over k}, ~~~~~~~|x|<{\pi \over 2},$$
wherein $B_n$ are the Bernoulli numbers.
For R10 and R12 one may apply, together with double angle formulas,
$$\sum_{k=0}^\infty (-p)^k \cos kx={{1+p\cos x} \over {1+2p\cos x+p^2}},$$
with $p\to \exp(a\cos^2\theta)$ and $x \to a\sin \theta \cos \theta$.
R10 reads
$$\int_0^{\pi/2} {{1+e^{a\cos^2\theta}\cos(a\sin \theta \cos\theta)} \over {1+2e^{a\cos^2 \theta}
\cos(a\sin \theta \cos \theta)+e^{2a\cos^2 \theta}}}d\theta={\pi \over {2(e^{a/2}+1)}}.$$
Due to the periodicity and symmetry of the integrand, $\int_0^{\pi/2}(\ldots)d\theta={1 \over 4}\int_0^{2\pi}(\ldots)d\theta$.
Very recently other integrals with trigonometric-hyperbolic function integrands have been
discussed in \cite{coffey2018}.


\medskip
{\bf 3.  Generalized Russell-style integrals}
\medskip

In referring to entries of the list of \cite{tewodros1}, the integral identities are labelled with
a preceding A.  We present several generalizations.  
The entries of \cite{tewodros1} are said to result from the application of an integral
transformation with an integrand with a Gaussian factor and a kernel function of either $1/\cosh y$
or $1/\sinh y$.  Instead, many of these identities may be obtained with
$$2\sum_{j=0}^\infty e^{-(2j+1)\kappa x}={1 \over {\sinh \kappa x}}, ~~~~~~
2\sum_{j=0}^\infty (-1)^j e^{-(2j+1)\kappa x}={1 \over {\cosh \kappa x}}$$
and the interchange of summation and integration.

A1 is generalized to
$$\int_0^\infty x\left({{g\sinh gx} \over {\cosh^2 gx}}e^{-x^2/\pi^2}+{{\sqrt{\pi}\sinh x} \over
{\cosh^2 x}}e^{-g^2x^2}\right)dx=\int_0^\infty {e^{-x^2/\pi^2} \over {\cosh gx}}dx.$$
In \cite{tewodros1}, this identity is stated only for the specific value of $g=\gamma$, the
celebrated Euler constant \cite{lagarias}.  

Similarly, A10 is generalized to
$$\int_0^\infty \left(\pi^5 e^{-\pi^3 x^2/g}+g^{5/2}e^{-gx^2/\pi}\right){{x^2 dx} \over {\cosh \pi x}}
={{\pi g^{3/2}} \over 2}\int_0^\infty {e^{-g x^2/\pi} \over {\cosh \pi x}}dx,$$
whereas in \cite{tewodros1} $g$ is restricted to the value of the famous Catalan constant
$G=\sum_{k=0}^\infty {{(-1)^k} \over {(2k+1)^2}}.$

A3 is generalized to
$$\int_0^\infty x\left(e^{-x^2/\pi}+re^{-r^4 x^2/\pi}\right){{\sinh rx dx} \over {\cosh^2 rx}}
=\int_0^\infty {{e^{-r^2 x^2/\pi} dx} \over {\cosh r^2 x}},$$
whereas in \cite{tewodros1} $r=2$.

A5 is generalized to
$$\int_0^\infty xe^{-x^2/\pi} {{\sinh x dx} \over {\cosh^2 x}}={r \over 2}\int_0^\infty {e^{-r^2 x^2/\pi}
\over {\cosh r x}}dx.$$
In \cite{tewodros1}, this identity is stated only for $r=2$.

In regard to A6 we remark that, by using the substitution $x \to 1/x$,
$$\int_0^1 {{x^{-\ln x} dx} \over {1+x^2}}=\int_1^\infty  {{x^{-\ln x} dx} \over {1+x^2}}.$$
We then give the generalization
$$\int_0^1 {{x^{-\ln x} dx} \over {1+x^2}}={r \over 2}\int_0^\infty {e^{-r^2 x^2/\pi} \over {\cosh
r \sqrt{\pi} x}}dx,$$
rather than simply $r=2$ as in \cite{tewodros1}.

A9 is generalized to
$$\int_0^\infty \left(\sqrt{\pi} e^{-x^2/p}+p^2\sqrt{p}\pi^{-2}e^{-px^2/\pi^2}\right){{x^2 dx} \over
{\cosh x}}={{p\pi \sqrt{p}} \over 2}\int_0^\infty {{e^{-px^2} dx} \over {\cosh \pi x}},$$
while in \cite{tewodros1} this identity is stated just for $p=3$.

\pagebreak
{\bf 4.  Further directions}
\medskip

There exists many other analogous integrations, of which we mention only a few other
approaches.  One is q-extensions, in some cases using the q-binomial series.  Since the $_2F_1$
function (and so many of its specializations) have appeared, this suggests in addition a consideration of
finite field analogs. 

Bessel functions provide an analog of trigonometric functions, and although oscillatory, are not
periodic.  A better analogy to many Russell-type integrals is provided by elliptic functions, of
which there are many choices.  This opens many avenues for investigation.

One example is to consider the set of Jacobi elliptic functions.  A general idea is to consider
definite integrations over an interval(s) including rational multiples of one of the two periods.
We hereby provide two very specialized analogs of R5, with further $x=1$.  Jacobi elliptic functions
are denoted $sn(z,m)$, etc.\ with modulus $m$.  

In preparation, we have the particular complete elliptic integral {\bf K} of the first kind constant
$${\bf K}\left({1 \over 2}\right)={{\Gamma^2(1/4)} \over {4 \sqrt{\pi}}}.$$
Then
$$\int_0^{{\bf K}(1/2)} z ~sn\left(z,{1 \over 2}\right)\sqrt{1+cn^2\left(z,{1 \over 2}\right)}dz={\pi \over 2},$$
and
$$\int_0^{2{\bf K}(1/2)} z ~sn\left(z,{1 \over 2}\right)\sqrt{1+cn^2\left(z,{1 \over 2}\right)}dz
={{\Gamma^2(1/4)} \over {\sqrt{2\pi}}}.$$
Other possibilities are to consider various analogs of Lemma 1 in anticipation of series expansions of
integrands of rational form in elliptic functions.  Two of many such families of definite integrals are to
accordingly consider 
$$\int_0^{2{\bf K}(m)}z ~sn(z,m)dn(z,m)cn^k(z,m) dz, ~~~~~~\int_0^{{\bf K}(m)}z ~sn(z,m)dn(z,m)cn^k(z,m) dz,$$
with again $2\bf{K}$ being a half period of $cn$ and $sn$, $sn^2 z+cn^2 z=1$, and ${d \over {dz}}cn z=-sn z ~dn z$.  By way of integration by parts, there are again various generalizations of Lemma 1, and further
integrals with integrands of higher powers of $sn$ in terms of recursion relations.

We present the following examples for $m=1/2$, being again the lemniscatic setting.
$$\int_0^{{\bf K}(1/2)}z ~sn(z,1/2)dn(z,1/2)cn^2(z,1/2) dz={1 \over {3\sqrt{2}}},$$
$$\int_0^{{\bf K}(1/2)}z ~sn(z,1/2)dn(z,1/2)cn^4(z,1/2) dz={\pi \over {20\sqrt{2}}},$$
and
$$\int_0^{{\bf K}(1/2)}z ~sn(z,1/2)dn(z,1/2)cn^5(z,1/2) dz={1 \over {10\sqrt{2}G}},$$
where $G \approx 0.83462684167...$ is the Gauss constant, given by
$$G={\sqrt{2} \over \pi}{\bf K}\left({1 \over \sqrt{2}}\right)={1 \over {(2\pi)^{3/2}}}\Gamma^2\left({1 \over
4}\right).$$




\bigskip
\centerline{\bf Acknowledgements}
T. Amdeberhan is thanked for useful correspondence, and P. Martin for reading the manuscript.

\pagebreak

\end{document}